\documentclass[12pt]{amsart}

\usepackage{amssymb,latexsym,amsmath, amscd}

\usepackage{graphics}
\usepackage{epsfig}

\newtheorem{theorem}{Theorem}
\newtheorem{corollary}[theorem]{Corollary}
\newtheorem{dfn}{Definition}[section]

\newtheorem{rem}[dfn]{Remark}
\newtheorem{thm}[dfn]{Theorem}

\newtheorem{lem}[dfn]{Lemma}

\def\proof{\par\medskip\noindent{\it Proof. }}

\def\Z{{\mathbb Z}}

\def\al{\alpha}

\def\De{\Delta}

\def\la{\lambda}

\def\8{\infty}
\def\<{\langle}
\def\>{\rangle}

\def\BE{\begin{equation}}
\def\EE{\end{equation}}

\def\8{\infty}

\def\3{\ss}

\def\la{\lambda}
\def\lam{\lambda}

\def\>{\rangle}
\def\<{\langle}

\def\3{\ss}
\def\8{\infty}

\def\mf{\mathfrak }

\begin{document}
\def\p{{\mf p}}

\title{Saturation and Irredundancy for Spin(8)}

\author{Michael Kapovich, Shrawan Kumar and John J. Millson}
\thanks{The authors were partially supported by the NSF individual
grants DMS 0405180, DMS 0401084 and DMS 0405606 respectively and together
by an NSF FRG grant}

\begin{abstract}
We explicitly calculate the {\em  triangle
inequalities} for the group $PSO(8)$,  thereby explicitly solving
the {\em eigenvalues of a sum} problem for this group (equivalently
describing the side-lengths of geodesic triangles in the corresponding symmetric space for the
metric $d_{\Delta}$ with values in the Weyl chamber $\Delta$). We then apply
some computer programs to verify  two basic
questions/conjectures. First, we verify that  the above system of inequalities
is irredundant. Then, we verify the ``saturation conjecture'' for  the
decomposition of tensor products of finite-dimensional irreducible
representations of  $Spin(8)$. Namely, we show
that for any  triple of dominant weights $(\la, \mu, \nu)$ such
that $\lam+\mu+\nu$ is in the root lattice, and any positive integer $N$,
$$
(V(\lam)\otimes V(\mu)\otimes V(\nu))^{Spin(8)}\ne 0$$
if and only if 
$$
(V(N\lam)\otimes V(N\mu)\otimes V(N\nu))^{Spin(8)}\ne 0.
$$
\end{abstract}

\dedicatory{Dedicated to F. Hirzebruch on the occasion of his seventieth birthday}

 \maketitle

\section{Introduction}
\label{introduction}

In this paper we address the following three basic problems in  algebraic group
theory. The statements of the first two problems and the description of their solution set, the cone
$\mathcal{C}(R)$, do not depend on the fundamental group of $G$.

Let $G$ be a connected  complex semisimple  algebraic group. We fix a Borel subgroup
$B$,  a maximal torus $T \subset B$ and a maximal compact subgroup $K$.
Let $X= G/K$ be the associated symmetric space.
Let $\mathfrak{b}, \mathfrak{h}$, $\mathfrak{k}$ and  $\mathfrak{g}$ be the
Lie algebras of $B,T,K$ and $G$ respectively.  Let $\mathfrak{g} = \mathfrak{k} \oplus
\mathfrak{p}$ be the Cartan decomposition. We can and will  assume that
$\mathfrak{h}$ satisfies
$\mathfrak{h} = \mathfrak{h}\cap \mathfrak{k} \oplus \mathfrak{h} \cap \mathfrak{p}.$
Let $\mathfrak{a}$ be the second intersection $\mathfrak{h} \cap \mathfrak{p}$
(the Cartan subspace). Let $A$ be the real split subtorus of $T$
corresponding to $\mathfrak{a}$.
The choice of  $B$ determines the set $R^+ \subset \mathfrak{a}^*$ of positive
roots and thus the set $\Pi =\{\alpha_1, \dots,
\alpha_l\}\subset R^+$ of simple roots and also the fundamental
weights $\{\omega_1, \dots, \omega_l\}, \,l$ being the rank of
$G$. The cone generated by the  positive roots
determines the dual cone $\Delta \subset \mathfrak{a}$, the (closed) Weyl chamber.

In Section \ref{background} we will introduce the $\Delta$-valued
distance $d_{\Delta}$ on the symmetric space $X$. We then have the following

\medskip
{\bf Problem 1.} {\em The generalized triangle inequalities}.
Give conditions on a triple $(h_1,h_2,h_3) \in \Delta^3$ that are
necessary and sufficient in order that there exist a triangle in
$X$ with vertices $x_1,x_2,x_3$ such that $d_{\Delta}(x_1,x_2) =
h_1$, $d_{\Delta}(x_2,x_3) = h_2$ and $d_{\Delta}(x_3,x_1) =
h_3$.

\medskip

Our second problem  is the  generalization (to
general $G$)
of the problem of finding the possible eigenvalues of a sum
of Hermitian matrices given  the eigenvalues of the summands.
To formulate this problem, define
the map
$$
\pi: {\mf p}/K \to \Delta
$$
 by taking $\pi(x)$ to be the unique point in the intersection of $\Delta$ with the
 $Ad\,K$-orbit of $x$.

\medskip
{\bf Problem 2.} {\em The generalized eigenvalues of a sum
problem}.
Determine the subset ${\mathcal C}={\mathcal C}(R)\subset
\Delta^3$ consisting of triples  $(h_1, h_2,
h_3)\in \Delta^3$ such that there exists a triple  $(y_1, y_2, y_3)\in {\mf
p}^3$ for which
$$
y_1+y_2+y_3=0
$$
and $\pi(y_1)=h_1, \pi(y_2)=h_2, \pi(y_3)=h_3$.

\medskip
It turns out that the sets of solutions to Problems 1 and 2 coincide, see
\cite{Klyachko2},\
\cite{AlekseevMeinrenkenWoodward},\ \cite{EvensLu},\ \cite{KapovichLeebMillson1}
and \cite{KapovichLeebMillson2}. The common solution set  ${\mathcal C}$ is in fact a convex
 homogeneous polyhedral cone ${\mathcal C}$, \cite{BerensteinSjamaar},
see also \cite{KapovichLeebMillson1}. The set
${\mathcal C}$ is described in \cite{BerensteinSjamaar} with a refinement in \cite{KapovichLeebMillson1} by a certain system
of homogeneous linear inequalities,
the {\em extended triangle inequalities} $\widetilde{TI}(R)$, which
is, in general, a redundant system.  A smaller system,  the {\em triangle
inequalities} $TI(R)$, was introduced by Belkale and Kumar in
\cite{BelkaleKumar}. These systems of inequalities are based on the cup
product, resp. the degenerated cup product, on the cohomology of the
generalized Grassmannians $G/P$, where $P\subset G$ are maximal
parabolic subgroups.

Belkale and Kumar have posed the question if the system $TI(R)$ is irredundant
(cf. \cite{BelkaleKumar}, Section 1.1).

\begin{rem}
The two systems of inequalities  $\widetilde{TI}(R)$ and ${TI}(R)$ coincide in the case of type $A$ root
systems (cf. \cite{BelkaleKumar}, Lemma 19). In this case irredundancy was proven by Knutson, Tao and
Woodward in \cite{KnutsonTaoWoodward}.
\end{rem}

 In this paper we prove that the system ${TI}(R)$ is indeed irredundant
for $R=D_4$, that is
for groups with Lie algebra $\mathfrak{so}(8)$.

\begin{thm} The system of triangle inequalities
$TI(D_4)$ is irredundant.
\end{thm}

We will prove this by explicitly computing the system $TI(D_4)$ and
then verifying that it consists of $306$ inequalities, while the
cone ${\mathcal C}$ has $306$ facets (and $81$ extremal rays).
This computation is done by applying the computer program CONVEX,
\cite{convex}.

\medskip
Our third problem concerns the decomposition of tensor products of finite-dimensional irreducible
representations of  a complex semisimple group. It is natural to assume that
this group is simply-connected.

In order to
relate the third problem to the first two it is necessary to introduce
the Langlands' dual $G^{\vee}$ of $G$, see \cite{Springer}, pages 3-6. We explain briefly why this is the
case.
There is a natural correspondence of maximal tori $T$ and $T^{\vee}$ for the two
groups such that the
 dominant coweights  (the ``integral points'' in $\Delta$) of $G$ are the
dominant weights  of $G^{\vee}$ whence {\em the input data for Problem 3
for the case of $G^{\vee}$ is a subset consisting of the ``integral points"  of
the input data for Problems 1 and 2 for  $G$}.

{\it For this reason we now assume that in the previous discussion $G$
was the centerless form of $\mathfrak{g}$ (since this corresponds to
the assumption that  the group $G^{\vee}$
is simply-connected).}

Let $P^{\vee}$ be the weight lattice of $G^{\vee}$, i.e.,
$P^{\vee}$ is the character lattice of the maximal torus $T^{\vee} \subset G^{\vee}$. Then,
$$
{\mathcal D}={\mathcal D} (G^{\vee}):=P^{\vee}\cap \Delta
$$
is the set of dominant weights of $G^{\vee}$.

\begin{dfn}
We define $(\mathcal{D}^3)^0$ to be the subsemigroup of $\mathcal{D}^3$
consisting of those triples of dominant weights whose sum is in
the root lattice $Q^\vee$  of $G^\vee$.
\end{dfn}

Given $\la\in {\mathcal D}$
let $V(\la)$ denote the irreducible representation of $G^\vee$ with
dominant weight $\la$.

We now state our third problem.

\medskip
{\bf Problem 3.} Determine the semigroup $\mathcal{R}=\mathcal{R}(G^\vee)\subset
{\mathcal D}^3$ consisting of triples of dominant weights $(\la,
\mu, \nu)$ such that
$$
(V(\la)\otimes V(\mu)\otimes V(\nu))^{G^{\vee}}\ne 0.
$$

It is well known that Problems 2 and 3 are related,
namely that the semigroup
$\mathcal{C}^0=\mathcal{C}(G^\vee)^0 :=\mathcal{C}(R)\cap ({\mathcal{D}^3})^0$
is the saturation of the semigroup
$\mathcal{R}(G^\vee)$ in the semigroup  $({\mathcal{D}^3})^0$.
For a more detailed statement,  see
Theorem \ref{problemsoneandthree}.

It was conjectured in \cite{KapovichMillson2} that ${\mathcal R}$
is saturated in $({\mathcal{D}^3})^0$ for all the simply--laced (and simply-connected)
groups $G^\vee$, i.e., for any $(\lam,\mu,\nu)\in ({\mathcal{D}^3})^0$ and any
positive integer $N$, if $(N\lam,N\mu,N\nu)\in {\mathcal{R}}$, then
$(\lam,\mu,\nu)\in {\mathcal{R}}$, i.e.,  $\mathcal{R}(G^\vee)=\mathcal{C}(G^\vee)^0$.
This conjecture is again
known in the case of type $A$ root systems, this is the {\em
saturation theorem} of Knutson and Tao \cite{KnutsonTao}, see also
\cite{Belkale, DerksenWeyman, KapovichMillson1} for alternative proofs.

We now state our second main theorem.

\begin{thm}\label{trianglerep}
Let $R=D_4$ so that $G^{\vee}=Spin(8)$. A triple  $(\lambda,\mu,\nu)\in
({\mathcal{D}^3})^0$ satisfies   $(\lambda,\mu,\nu)\in {\mathcal{R}}(Spin (8))$ if and only if
 $(\lambda,\mu,\nu)\in {\mathcal{C}}(D_4)$.
Equivalently,  the semigroup $\mathcal{R}(Spin (8))$ is saturated in
the semigroup $({\mathcal{D}^3})^0$.
\end{thm}

In order to prove Theorem \ref{trianglerep} we use the computer program 4ti2,
\cite{42}, to compute the Hilbert basis of the semigroup
${\mathcal C}(Spin (8))^0$. It turns out that this basis consists of 82
elements (just one more than the number of extremal rays). Moreover,
modulo the permutations of the vectors $\la, \mu, \nu$ and action of the
automorphisms of the Dynkin diagram of $D_4$, there are only $10$
different semigroup generators. For each of these generator $(\la_i,
\mu_i, \nu_i)$ we verify that
$$
(\la_i, \mu_i, \nu_i)\in {\mathcal R}
$$
by applying the MAPLE package WEYL, \cite{stembridge}. Since ${\mathcal
R}$ is a semigroup, it then follows that ${\mathcal
C}^0={\mathcal R}$.

By \cite{KapovichLeebMillson3}, Section 9.4, the previous theorem implies
the following  saturation
theorem for the structure constants of the spherical Hecke algebra
of $PSO(8)$. Considering $PSO(8)$ as a group scheme  $\underline{PSO(8)}$ over $\mathbb{Z}$,
let ${\mathcal G}$ be the group of  its rational points in a nonarchimedean local field $\mathbb{K}$. Let $\mathcal{O}$ be the
ring of integers (elements of nonnegative valuation) of $\mathbb{K}$. We  let ${\mathcal K}$
be the group of $\mathcal{O}$-rational points of $PSO(8)$. Let $\mathcal{H}_{\mathcal G}$
denote the associated spherical Hecke ring. We recall that the set
of dominant coweights $\mathcal{D}$ of  $G$ parametrizes the ${\mathcal K}$-double
cosets in ${\mathcal G}$ and that the ring $\mathcal{H}_{\mathcal G}$ is free over $\mathbb{Z}$
with basis the characteristic functions $\{f_{\lambda}:\lambda \in \mathcal{D} \}$.
We let $*$ denote the (convolution) product in $\mathcal{H}_{\mathcal G}$.
We have

\begin{thm}\label{triangleHecke}
Let $G=PSO(8)$. For $\lambda,\mu,\nu \in {\mathcal D}$, the characteristic function of the
identity ${\mathcal K}$-double coset occurs in the expansion of  the product $f_{\lambda} * f_{\mu} *f_{\nu}$
in terms of the above basis
if and only if the triple $(\lambda,\mu,\nu)\in {\mathcal C}$
and $\lambda + \mu + \nu$ is in the coroot lattice $Q^\vee$ of $G$.
\end{thm}

\section{Further discussion of the three problems}
\label{background}

In this section we give some more details about the three problems formulated in the
Introduction. We follow the same notation (as in the Introduction). In particular,
$G$ is a complex semisimple adjoint group (with root system $R$) and $G^\vee$ is its
Langlands' dual, which is simply-connected (since $G$ is adjoint).

\subsection{The distance $d_{\Delta}$}

 We now define the $\Delta$-valued distance $d_{\Delta}$.
Let $A_{\Delta} $ be the
image of $\Delta$ under the exponential map $\exp:\mathfrak{g} \to
G$.   We will need the following basic
theorem, the Cartan decomposition for the group $G$, see
\cite{Helgason}, Theorem 1.1, pg. 402.

\begin{thm}\label{Cartandecomposition}
We have
$$G = KA_{\Delta}K.$$
Moreover, for any $g \in G$, the intersection of the double coset
$KgK$ with $A_{\Delta}$ consists of a single point to be denoted
$a(g)$.
\end{thm}

Let $\overline{x_1x_2}$ be the oriented geodesic segment in
$X=G/K$ joining the point $x_1$ to the point $x_2$. Then there exists
an element $g \in G$ which sends $x_1$ to the base point $o=eK$ and $x_2$ to
$y = \exp(\delta)$ where $\delta \in \Delta$. Note that the point
$\delta$ is uniquely determined by $\overline{x_1x_2}$. We define
a map $\sigma$ from $G$-orbits of oriented geodesic segments to
$\Delta$ by

$$\sigma(\overline{x_1x_2}) = \delta.$$

Clearly we have the following consequence of the Cartan
decomposition.

\begin{lem}
The map $\sigma$ gives rise to a one-to-one correspondence between
the $G$--orbits of oriented geodesic segments in $X$ and the points of
$\Delta$.
\end{lem}

In the rank $1$ case $\sigma(\overline{x_1x_2})$ is just
the length of the geodesic segment $\overline{x_1x_2}$. 

\begin{dfn}
We  call $\sigma(\overline{x_1x_2})$ the $\Delta$--length of
$\overline{x_1x_2}$ or the $\Delta$--distance between $x_1$ and
$x_2$. We  write $d_{\Delta}(x_1,x_2) =
\sigma(\overline{x_1x_2})$.
\end{dfn}

 We note the formula
 $$d_{\Delta}(x_1,x_2) = \log a(g_1^{-1}g_2) \ \text{where} \ x_1 = g_1K,
x_2 =g_2K.$$

\begin{rem}
The delta distance is {\em symmetric} in the sense that
$$d_{\Delta}(x_1,x_2) = - w_o d_{\Delta}(x_2,x_1),$$
where $w_o$ is the unique longest element in the  Weyl group.
\end{rem}

\subsection{The relations between Problems 1, 2 and 3}
In this subsection we expand the discussion in the Introduction concerning
the relations between the three problems. We first discuss the
relation between Problems 2 and 3. 

The (a)-part of the following theorem is standard, see for example the appendix of
\cite{KapovichLeebMillson3}. The (b)-part follows from Theorem 1.2 of \cite{KapovichLeebMillson1}, see also  Theorem 1.3 of \cite{KapovichLeebMillson2} and the paragraph following it. Of course,  the (b)-part is clear for the simply-laced groups. So, the only nontrivial case  is essentially that of the group $G$ corresponding to the root systems of type $B_l$. In this case, Belkale and Kumar have shown that the triangle inequalities themselves coincide under the identification of ${\mathfrak a}$ with ${\mathfrak a}^*$ (via any Killing form). 

\begin{thm}\label{problemsoneandthree}
(a) For any semisimple adjoint group $G$ with root system $R$, 
under the identification of ${\mathfrak a}$ with ${\mathfrak a}^*$ (via any Killing form), 
 $$ {\mathcal R}(G^\vee)\subset {\mathcal C}(R^\vee).$$
 Conversely, for any triple $(\lambda,\mu,\nu)\in  {\mathcal C}(R^\vee)\cap
 {\mathcal D}^3$,
there exists a positive integer $N$ such that $(N\lambda, N\mu, N\nu) \in {\mathcal R}(G^\vee).$

(b) Under the identification of ${\mathfrak a}^*$ with ${\mathfrak a}$,
 $$ {\mathcal C}(R)= {\mathcal C}(R^\vee).$$

Thus, combining the two parts, we get the following  intrinsic inclusion:
 $$ {\mathcal R}(G^\vee)\subset {\mathcal C}(R).$$
\end{thm}

We recall the following standard definition.  
Suppose that $S_1 \subset S_2$ is an inclusion of semigroups. Then
the {\em saturation} of $S_1$ in $S_2$ is the semigroup of elements
$x \in S_2$ such that there exists $n \in \mathbb{Z}_+$ with
$nx \in S_1$.

\begin{figure}[tbh]
\centerline{\epsfxsize=3.5in \epsfbox{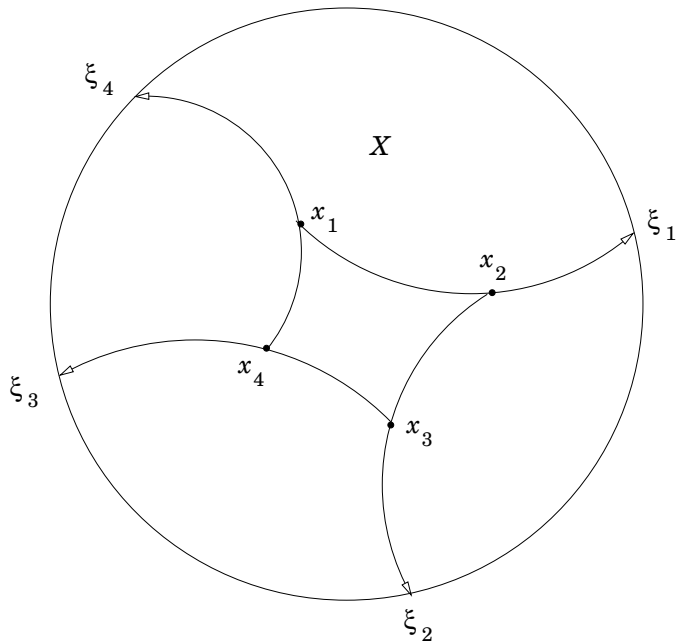}}
\caption{\sl Gauss map.}
\label{gauss.fig}
\end{figure}

\begin{rem}
We may restate the previous theorem by saying that the semigroup
$\mathcal{C}(G^\vee)^0:={\mathcal C}(R)\cap ({\mathcal D}^3)^0$ is the saturation of the semigroup
$\mathcal{R}(G^\vee)$ in the semigroup  $({\mathcal{D}^3})^0$.
\end{rem}

We conclude this section by briefly indicating why the solutions of
Problems 1 and 2 coincide. First of all, Problem 2 (in the case of an $n$--fold sum) 
can be reformulated geometrically as a problem of the existence of geodesic polygons in $\p$ with fixed 
$\De$--valued side-lengths. Both $\p$ and the symmetric space $X$ admit compactifications by a ``visual'' \
sphere $S$ which also has the structure of a spherical building 
$\partial_{\text{Tits}} X$. The vertices of 
this building are points in the flag manifolds $G/P$ (where $P$'s are the maximal parabolic subgroups in $G$).

Then one can define the {\em Gauss map} $\Gamma$ which sends the geodesic polygon  
$[x_1,...,x_n]\subset X$ to the weighted configuration 
$$
\Gamma([x_1,... , x_n])=((m_1, \xi_1),... , (m_n, \xi_n))
$$
of points in $S$. Here $m_i:=d(x_i, x_{i+1})$ are the ordinary distances, which serve as {\em weights} at the points 
$\xi_i\in S$.  The same definition also works for $X$ replaced by  
 $\p$. The Gauss map from quadrilaterals in the hyperbolic plane to configurations of four points $\xi_1,\xi_2,
\xi_3,\xi_4$ on the visual boundary (the circle) is depicted in figure \ref{gauss.fig}. The key problem 
then is to identify the images of Gauss maps $\Gamma$ for $X$ and  $\p$. It turns out that both 
consist of nice   semistable weighted configurations on $S$, where the notion of stability 
is essentially the one introduced by Mumford in Geometric Invariant Theory 
(in the case when the weights are natural numbers).  
We refer the reader to Theorems 5.2 and 5.9 of \cite{KapovichLeebMillson1} for the precise statements. 

Therefore, Problems 1 and 2 are both equivalent to the existence problem for semistable weighted configurations 
on $S$ and hence Problems 1 and 2 are equivalent.

\section{The triangle inequalities}

 We  need more
notation. We let $S=\{s_1,\dots, s_l\}$ be the set of (simple) reflections in the root
hyperplanes defined by the simple roots and let $W\subset
\text{Aut}\, \mathfrak a$ be the Weyl group generated by $S$.

Let $\{x_i\}$ be the basis of $\mathfrak{h}$ dual to the basis  $\Pi$,
i.e., $\alpha_i(x_j)=\delta_{i,j}$.  We let $\ell$ be
the length function on $W$. Let
$\alpha_i^\vee$ be the coroot corresponding to the root $\alpha_i$.
Also, for  a standard parabolic subgroup $P$ of $G$ (i.e. $P\supset B$), we let $W_P$ be the subgroup
of elements with representatives in $P$ and $W^P$ denote
the set of shortest length representatives for the cosets $W/W_P$
(we recall that each coset has a unique shortest length representative).
Let $w^P_o$ be the unique longest element in $W^P$.

\subsection{The extended triangle inequalities}

We now describe the solution of Problem 1 of the Introduction, that
is the description of the inequalities determining the ${\Delta}$-valued
side-lengths of geodesic triangles in $X$.

\subsubsection{The weak triangle inequalities}

We first describe a natural subsystem
of the triangle inequalities.
The naive triangle inequality
$$d_{\Delta}(x_1,x_3) \leq_{\Delta} d_{\Delta}(x_1,x_2) + d_{\Delta}(x_2,x_3)$$
does not hold \cite{KapovichLeebMillson3}. Here the order $\leq_{\Delta}$ is the
one defined by the (acute) cone $\Delta$. This can be remedied  if we replace $\Delta$
by the dual (obtuse) cone $\Delta^*$ and let $\leq_{\Delta^*}$
denote the associated order.
Then, the analogue of the above inequality holds and, in fact,
for any element $w$  of the  Weyl group $W$, the inequality
$$w \cdot d_{\Delta}(x_1,x_3) \leq_{\Delta^*} w \cdot d_{\Delta}(x_1,x_2) + d_{\Delta}(x_2,x_3)$$
holds.
We call the resulting system of inequalities (as $w$ varies)
the {\em weak triangle inequalities} to be denoted $WTI(R)$.

For the root systems $R$ of ranks one and two, the weak triangle inequalities  already
give a solution to
Problems  1 and 2 of the Introduction. However,  they are no longer  sufficient in ranks three
or more.

\subsubsection{The extended triangle inequalities}

We now describe  a  system of linear inequalities on $\mathfrak{a}^3$
which describes the cone ${\mathcal C}(R)$. However, this system is
usually not  irredundant. These inequalities (based on the cup-product of Schubert classes) will be called the {\em extended triangle
inequalities}. The system of extended triangle inequalities is independent of the choice of $G$
corresponding to a fixed Lie algebra $\mathfrak{g}$,  hence  depends only on
the root system $R$ associated to $G$.  We  denote the system of
extended triangle inequalities by $\widetilde{TI}(R)$ or just $\widetilde{TI}$
when the reference to $R$ is clear.

As a consequence of  the Bruhat decomposition,
\[G =  \bigsqcup_{w \in W^P} BwP,\]
 the generalized flag variety $G/P$ is
the disjoint union of the subsets 
$$
\{C_w^P:= BwP/P\}_{w\in W^P}.
$$
The subset $C_w^P$ is biregular isomorphic  to the affine space
${\mathbb C}^{\ell(w)}$ and is called a {\em Schubert cell}, where $\ell(w)$
is the length of $w$. The
closure $X^P_w$ of $C_w^P$ is called a {\em Schubert variety}. We
will use $[X_w^P]$ to denote the integral homology class in
$H_*(G/P)$ carried by $X_w^P$.
Then,
the integral homology $H_*(G/P)$ is a free $\Z$--module with basis
$\{[X_w^P]: w \in W^P\}$.

Let  $\{
\epsilon_w^P : w \in W^P\}$  denote the dual basis of $H^*(G/P)$ under the Kronecker
pairing $\<\ ,\ \>$ between homology and cohomology.
Thus, we have
for $w,w^{\prime}\in W^P,$

$$\<\epsilon_w^P, [X^P_{{w^{\prime}}}]\> = \delta_{w,w^{\prime}}.$$

The system of extended triangle inequalities  breaks up into
rank$(\mathfrak{g})$ subsystems ${\widetilde{TI}}^P$, where $P$ runs over standard
maximal parabolic subgroups. The subsystem $\widetilde{TI}^P$ is controlled by
the Schubert calculus in the generalized Grassmannian $G/P$ in the
sense that there is one inequality $T^P_{\mathbf{w}}$ for each
triple of elements $\mathbf{w} =(w_1,w_2,w_3) \in W^P$ such that
$$\epsilon^P_{w_1}\cdot \epsilon^P_{w_2} \cdot \epsilon^P_{w_3} = \epsilon^P_{w_o^P}$$
in  $H^*(G/P)$.  To describe the inequality
$T^P_{\mathbf{w}}$,  let $\lambda$ be the  fundamental weight corresponding to $P$. Then the action of
$W$ on $\mathfrak{a}^*$ induces a one-to-one
correspondence $f:W^P \to W\lambda$. Thus,  we may reparameterize the
Schubert classes in $G/P$ by elements of $W\lambda\subset \mathfrak{a}^*$. We let
$\lambda_i = f(w_i), i = 1,2,3$.  Then the inequality
$T^P_{\mathbf{w}}$ is given by
$$\lambda_1(h_1) + \lambda_2(h_2) + \lambda_3(h_3) \geq 0,
\ (h_1,h_2,h_3) \in \Delta^3.$$

\subsection{The triangle inequalities}

As we have mentioned earlier, the system of  extended  triangle inequalities is in general
not an irredundant system. We now describe the subsystem of triangle inequalities.

To this end
we recall the definition of the  new product
$\odot_0$ in the cohomology $H^*(G/P)$ introduced by Belkale-Kumar [BK, Sect. 6].
We only need to consider the case when  $P$ is  a standard
maximal parabolic subgroup. In this case, we set $x_P=x_{i_P}$, where $s_{i_P}$ is
the unique simple reflection not in $W_P$. Similarly, we set $\omega_P=\omega_{i_P}.$
We can identify $W^P$ with the orbit
$W\cdot \omega_P$. For $w\in W^P$, let $\la_w=\la_w^P$ denote $w( \omega_P)$; this is
called the
{\it maximally singular weight} corresponding to $w$.

Write the cup product in $H^*(G/P)$ as follows:
\[\epsilon_u^P\cdot \epsilon_v^P=\sum_{w\in W^P}\,d^w_{u,v} \epsilon_w^P.\]
Then, by definition,
\[\epsilon_u^P\odot_0 \epsilon_v^P=\sum_{w\in W^P}\,d^w_{u,v}
\delta^w_{u,v}\epsilon_w^P,\]
 where $\delta^w_{u,v}:=1$ if $
(u^{-1}\rho+ v^{-1}\rho-
w^{-1}\rho-
\rho)(x_P)=0$ and
 $\delta^w_{u,v}:=0$ otherwise, where $\rho$ is the (standard) half sum of positive roots
of $\mathfrak{g}$.

Recall that $\pi:\mathfrak{p}/K \to \Delta$
is defined by intersecting an $Ad K$-orbit with $\Delta$. Then [BK,
Theorem 28] gives the following solution of Problem 2 stated in the Introduction:   

\begin{thm}\label{eigen} Let $(h_1,\dots,h_n) \in\Delta^n$.
Then, the following are equivalent:

(a) There exists $(y_1,\dots,y_n)\in \mathfrak p^n$ such that $\sum_{j=1}^n
y_j=0\,\,\text{and }\, \pi(y_j)=h_j$ for all $j=1,\dots,n.$

(b)
 For every standard maximal parabolic subgroup $P$ in $G$ and every choice of
  $n$-tuple  ${\bf w}=(w_1, \dots, w_n)\in (W^P)^n$ such that
$$\epsilon^P_{w_1}\odot_0\, \cdots \,\odot_0 \epsilon^P_{w_n}=
\epsilon^P_{w_o^P} \in \bigl(H^*(G/P), \odot_0\bigr),$$
the following inequality holds:
\begin{equation*}\label{eqnB}
\sum_{j=1}^n\,\la^P_{w_j} (h_j)\geq 0.\tag{$T^P_{{\bf w}}$}
\end{equation*}
\end{thm}

The collection of inequalities $\{T^P_{{\bf w}}\}$, such that ${\bf w}$ and $P$ are  as in (b), is called the {\em triangle inequalities}.
 
\begin{rem}
As was the case for $n=3$, the statement in (a) is equivalent to
the existence of a geodesic $n$-gon in $X$ with $d_{\Delta}$--side-lengths
$h_1,h_2,\dots,h_n$.
\end{rem}

\section{Determination of the product $\odot_0$ in $H^*(G/P)$}

  \bigskip

{\it From now on,  the group $G$ will be taken to be the adjoint
group of type $D_4$, i.e.,} $G=PSO(8)$. Since $G$ is simply-laced, the Langlands' dual
$G^\vee$ has the same root system as $G$. Moreover, $G$ being the adjoint group,
$G^\vee$ is the simply-connected cover of $G$, i.e., $G^\vee = Spin (8)$. We will only
need to consider the maximal parabolic subgroups.  We will abbreviate the
classes $\epsilon_w^P$ for $w \in W^{P}$ by $b_i^j$ according to the
following tables. Here the subscript $i$ denotes  half of the
cohomological degree of $b_i^j$, i.e., $b_i^j \in
H^{2i}(G/P)$, and $j$ runs over the indexing set with cardinality
equal to the rank of $H^{2i}(G/P)$. In the case that $H^{2i}(G/P)$
is of rank one, we suppress the superscript $j$. Moreover, in the
following tables, we also list the  maximally singular weight
$\lambda_w:=w\omega_P$ associated to the element $w\in W^{P}$ as well as 
the value $n_w:= (w^{-1}\rho)(x_P)$. We express $\lambda_w$ in
terms of the standard coordinates $\{\epsilon_i\}_{i=1, \dots, 4}$
of $\mathfrak{h}^*$ as given in [Bo, Planche IV]. We follow the
following indexing convention as in loc cit.


\begin{figure}[tbh]
\centerline{\epsfxsize=3.5in \epsfbox{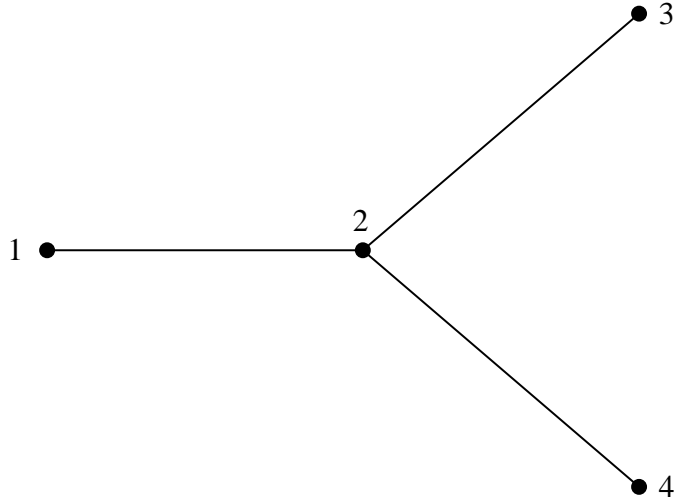}}
\caption{\sl Dynkin diagram for $D_4$.}
\label{Dynkin.fig}
\end{figure}

For any $1\leq i\leq 4$, let $P_i$ be the maximal parabolic subgroup of
$G$ corresponding to
the $i$-th node of the Dynkin diagram, i.e., $W_{P_i}$ is generated by the simple
reflections $\{s_j; j\neq i\}$.

To determine the cohomology $H^*(G/P_i)$ under the product $\odot_0$,
for all the maximal parabolic subgroups $P_i$,
because of the Dynkin automorphisms, we only need to determine it for $i=1,2$. Moreover,
since $P_1$ is a minuscule parabolic in Spin$(8)$, the product $\odot_0$
coincides with the cup product by [BK, Lemma 19].

In what follows we will use the symbol $F$ to denote the group
of automorphisms of the Dynkin diagram of $D_4$, so $F$ is
isomorphic to the symmetric group $S_3$ on the labels $\{1,3,4\}$ of the Dynkin diagram.

\subsection{Determination of $(H^*(G/P_1), \odot_0)$}

The longest element $w_o$ of $W$ is given by
\begin{equation}\label{eqn1}
w_o=s_4s_2s_1s_4s_2s_4s_3s_2s_4s_1s_2s_3,
\end{equation}
and it is central in $W$.
Moreover, the longest element $w_{o,P_1}$ of $W_{P_1}$ is given by
\begin{equation}\label{eqn2}
w_{o,P_1}=s_3s_2s_4s_3s_2s_3.
\end{equation}
Thus, the longest element $w_o^{P_1}$ of $W^{P_1}$ is given by (cf. [KuLM, Proposition
2.6])
\begin{equation}\label{eqn3}
w_o^{P_1}=w_ow_{o,P_1}=s_1s_2s_3s_4s_2s_1.
\end{equation}
From this and the fact that $\mid W^{P_1}\mid =8$, we see that the elements
of $W^{P_1}$ are enumerated as in the chart below. To calculate $n_w$, use the general formula
(cf. [K, Corollary 1.3.22]) for any $w\in W$:
\begin{equation}\label{eqn4}
\rho - w^{-1}\rho=\sum_{R^+\cap w^{-1}R^-}\,\alpha,
\end{equation}
and for any parabolic subgroup $P$ of $G$ and any $w\in W^P$,
\begin{equation}\label{eqn5}
R^+\cap w^{-1}R^-\subset R^+\setminus R^+_P,
\end{equation}
where $R^-:=-R^+$ and $R^+_P$ is the set of positive roots in the Levi component of $P$.
Since $P_1$ is a minuscule maximal parabolic subgroup, for any $w\in W^{P_1}$, by
(\ref{eqn4}) and (\ref{eqn5})  we get
\begin{equation}\label{eqn6}
(\rho-w^{-1}\rho)(x_{P_1})=\ell(w).
\end{equation}
From this, the value of $n_w$ given in the following  chart can easily be verified since
$\rho(x_{P_1})=3.$ The value of $\lambda_w$ is obtained by explicit calculations.

\begin{equation*}
\begin{array}{|c|c|c|c|}
\hline
\epsilon_w^{P_1} &  ~w& \la_w & n_w\\
\hline
b_0=1   &  e    & (1,0,0,0) & 3\\
\hline
\hline
b_1     &s_1  & (0,1,0,0) & 2\\
\hline
\hline
b_2     &s_2s_1  & (0,0, 1,0) & 1\\
\hline
\hline
b_3^1   &s_3s_2s_1  & (0,0, 0,1) & 0\\
\hline
b_3^2   &s_4s_2s_1  & (0,0, 0,-1) & 0\\
\hline
\hline
b_4     &s_3s_4s_2s_1 & (0,0, -1,0) & -1\\
\hline
\hline
b_5     &s_2s_3s_4s_2s_1   & (0,-1, 0,0) & -2\\
\hline
\hline
b_6     &s_1s_2s_3s_4s_2s_1   & (-1, 0, 0,0) & -3\\
\hline
\end{array}
\end{equation*}
\vskip3ex

Using [KuLM, Lemma 2.9] and the Chevalley formula (cf. [K, Theorem 11.1.7(i)]), all the
products in the following table can be determined except the products of $b_2$ with $b_2$ and
$b_3^*$. Since $b_1b_1=b_2$, using the Chevalley formula twice, we get  these products
as well.

\vskip2ex

{\bf Multiplication table for $G/P_{1}$ under the product $\odot_0$:}

\begin{equation*}
\begin{array}{|c||c||c||c|c|}
\hline
\odot_0     & b_1 & b_2 & b_3^1 & b_3^2 \\
\hline
\hline
b_1     &b_2        &       &  & \\
\hline
\hline
b_2     &b_3^1+ b_3^2   &  2b_4 & & \\
\hline
\hline
b_3^1   &b_4        & b_5   & 0  &\\
\hline
b_3^2   &b_4        &  b_5  &     b_6 & 0\\
\hline
\hline
b_4     &b_5 & b_6  &0 &0 \\
\hline
\hline
b_5     &b_6   &  0 & 0& 0\\
\hline
\hline
b_6     & 0 & 0 & 0& 0\\
\hline
\end{array}
\end{equation*}

\subsection{Determination of $(H^*(G/P_2), \odot_0)$}

For any parabolic subgroup $P$, let $\theta^P$ be the involution of $W^P$
defined by
\begin{equation}\label{eqn7}
\theta^Pw=w_oww_{o,P}.
\end{equation}
Then, by [KuLM, Section 2.1], $\epsilon^P_w$ is Poincar\'e dual to
$\epsilon^P_{\theta^Pw}.$

 Using
(\ref{eqn1})
and \begin{equation}\label{eqn8} w_{o,P_2}=s_1s_3s_4, \end{equation} we get
\begin{equation}\label{eqn9} w_o^{P_2}=s_2s_4s_1s_2s_3s_2s_1s_4s_2. \end{equation} The
enumeration of $W^{P_2}$ as in the following table can be read off from (\ref{eqn9}) together
with the fact that $\mid W^{P_2}\mid=24.$ The values of $\lambda_w$ and $n_w$ are
obtained by explicit calculations. Observe that the following identities provide some
simplification in the calculations of
 $\lambda_w$ and $n_w$: 

For any $w\in W^{P_2}$,
\begin{equation}\label{eqn10}
\lambda_{\theta^P(w)}=-\lambda_w,
\end{equation}
and
\begin{equation}\label{eqn11}
n_{\theta^P(w)}=-n_w.
\end{equation}

In the following table, the two $w$'s appearing in the same row are $\theta^P$--images of
each other, i.e., the corresponding classes $\epsilon^{P_2}_w$ are Poincar\'e dual to
each other.

{\scriptsize
 \begin{equation*}
\begin{array}{|c|c|c|c|c|c|c|c|}
\hline
\epsilon_w^{P_2} & ~w& \la_w & n_w &\epsilon_w^{P_2} & ~w& \la_w & n_w\\
\hline
b_0=1   &  e    & (1,1,0,0) & 5& b_9    &s_2s_4s_1s_2s_3s_2 s_1s_4s_2 &
(-1,-1,0,0) & -5\\
\hline
\hline
b_1     &s_2  & (1,0, 1,0) & 4&b_8  &s_4s_1s_2s_3s_2 s_1s_4s_2  &
(-1,0, -1,0) & -4\\
\hline
\hline
b_2^1   &s_1s_2  & (0,1, 1,0) & 3&b_7^1     &s_4s_2s_3s_2 s_1s_4s_2  &
(0,-1, -1,0) & -3\\
\hline
b_2^2   &s_3s_2  & (1,0, 0,1) & 3&b_7^2     &s_4s_2s_1s_2 s_3s_4s_2  &
(-1,0, 0, -1) & -3\\
\hline
b_2^3   &s_4s_2  & (1,0, 0,-1) & 3&b_7^3    &s_3s_2s_1s_2s_4s_3s_2
 & (-1,0, 0,1) & -3\\ \hline \hline b_3^1 &s_3s_1s_2 & (0,1, 0,1) & 2&b_6^1 &s_4s_2s_3s_1
s_4s_2 & (0,-1, 0,-1) & -2\\ \hline b_3^2 &s_4s_1s_2 & (0,1, 0,-1) & 2&b_6^2 &s_3s_2
s_4s_1s_3s_2 & (0,-1, 0,1) & -2\\ \hline b_3^3 &s_4s_3s_2 & (1,0, -1,0) & 2&b_6^3
&s_1s_2s_3s_4s_1s_2 & (-1,0, 1,0) & -2\\ \hline \hline b_4^1 &s_2s_3s_1s_2 & (0,0, 1,1) &
1&b_5^1 &s_4s_2s_1s_3s_2 & (0,0, -1,-1) & -1\\ \hline b_4^2 &s_4s_3s_1s_2 & (0,1, -1,0) &
1&b_5^2 &s_2s_4s_1s_3s_2 & (0,-1, 1,0) & -1\\ \hline b_4^3 &s_2s_4s_1s_2 & (0,0, 1,-1) &
1&b_5^3 &s_3s_2 s_1s_4s_2 & (0,0, -1,1) & -1\\ \hline b_4^4 &s_2s_4s_3s_2 & (1,-1,0,0) &
1&b_5^4 &s_1s_2s_3s_4s_2 & (-1,1, 0,0) & -1 \\ \hline \end{array} \end{equation*}
}

From the definition of $\odot_0$ and the values of $n_w$, we get the following.

  \begin{corollary}  \label{cor1} For $u,v \in W^{P_2}$, in $(H^*(G/P_2), \odot_0),$
\begin{align*}
\epsilon^{P_2}_u \odot_0 \epsilon_v^{P_2} &=
\epsilon^{P_2}_u \cdot \epsilon_v^{P_2}, \,\text{if \,}\,\ell (u)+\ell (v)
\leq 4\\
&=
\epsilon^{P_2}_u \cdot \epsilon_v^{P_2}, \,\text{if \,}\,\ell (u)+\ell (v)
\geq 5\,\text{and one of}\,\, \ell (u) \,\text{or}\,\, \ell (v)\geq 5
\\
&=
0, \,\text{if \,}\,\ell (u)+\ell (v)
\geq 5\,\text{and both of}\,\, \ell (u) \,\text{and}\,\,\ell (v)\leq 4.
  \end{align*}
  \end{corollary}

For any $i\neq j\in \{1,3,4\}$, let $\sigma_{i,j}$ be the involution of $H^*(G/P_2)$ induced from
the Dynkin diagram involution taking the $i$-th node to the $j$-th node and fixing
the other two nodes. Let $\hat{F}$ be the group of automorphisms  of $H^*(G/P_2)$  generated by
$\sigma_{1,3}, \sigma_{1,4}$ and $\sigma_{3,4}$. Then, $\hat{F}$ is isomorphic with the
symmetric group $S_3$.

Using [KuLM, Lemma 2.9], the Chevalley formula and Corollary
\ref{cor1}, we only need to calculate
$b_2^\cdot b_2^*, b_2^\cdot b_5^*, b_2^\cdot b_6^*$ and
$b_3^\cdot b_5^*.$ Further, using the automorphism group $\hat{F}$, it suffices to calculate
 $b_2^1 b_2^*, b_2^1 b_5^*, b_2^1 b_6^*$ and
$b_3^1 b_5^*.$ To calculate $ b_2^1 b_6^*$, write
$$ b_2^1 b_6^*=db_8, \,\,\,\text{for some}\, d.$$
Multiply this equation by $b_1$ and use the known part of the multiplication table to determine $d$. The calculation of
$b_3^1 b_5^*$ is exactly similar.

To calculate  $b_2^1 b_5^2$, write
$$ b_2^1 b_5^2=\sum_{i=1}^3 \,
d_ib_7^i, \,\,\,\text{for some}\, d_i\in\Z_+.$$
Multiplying the above equation by $b_1$, we get
 $$ b_2^1 b_5^2b_1=\sum_{i=1}^3 \,
d_ib_8.$$
On the other hand,
 $$ b_2^1 b_5^2b_1=\sum_{i=1}^3 \,
b_2^1b_6^i=2b_8.$$
Thus, $d_1+d_2+d_3=2.$ Using the involution $\sigma_{3,4}$ of $H^*(G/P_2)$, we are forced to
have
$$ b_2^1 b_5^2=b_7^2+b_7^3.$$
The calculation for  $b_2^1 b_5^4$ is similar. To calculate $b_2^1b_5^1$, write
$$
b_1b_1b_5^1=\sum_{i=1}^3\,b_2^ib_5^1.$$
 But,
$$b_2^3b_5^1=\sigma_{1,4}(b_2^1b_5^4)=0,$$
and
$b_2^2b_5^1=\sigma_{1,3}(b_2^1b_5^1).$ On the other hand
$$b_1b_1b_5^1=b_1b_6^1=b_7^1+b_7^2.$$
Thus,
$$ b_2^1 b_5^1=b_7^1 \,\,\text{or}\, b_7^2.$$
If $ b_2^1 b_5^1=b_7^2,$ then
$$b_2^1b_5^3=\sigma_{3,4}(b_2^1b_5^1)=b_7^3.$$
From this we conclude that $b^1_2b^1_2=0.$ However, by considering
the morphism $P_4/B \to \text{Spin}(8)/B$, induced from the
inclusion, we can easily see that   $b^1_2b^1_2\neq 0.$ This
contradiction forces $ b_2^1 b_5^1=b_7^1.$
 Using $\sigma_{3,4}$ as above,
we can calculate
 $b_2^1 b_5^3$ from  $b_2^1 b_5^1$.

To calculate  $b_2^1 b_2^*$, write
$$ b_2^1 b_2^*=\sum_{i=1}^4 \,
d_ib_4^i, \,\,\,\text{for some}\, d_i\in\Z_+.$$
Multiply this equation by $b_5^i$ to get
$$ b_2^1 b_2^*b_5^i=d_i.$$
Now, from the known part of the multiplication table, the $d_i$
can be determined.

\medskip
{\bf Multiplication table for $G/P_{2}$ under the product $\odot_0$:}

{\scriptsize
\begin{equation*}
\begin{array}{|c||c||c|c|c||c|c|c||c|c|c|c|}
\hline
\odot_0 & b_1 & b_2^1 & b_2^2 & b_2^3 & b_3^1 & b_3^2 & b_3^3 & b_4^1 & b_4^2 & b_4^3 &
b_4^4 \\
\hline
\hline
b_1 &  b_2^1+ b_2^2+ b_2^3 & & & & & & & & & & \\
\hline
\hline
b_2^1 &  b_3^1 + b_3^2 &  b_4^1 + b_4^3 & & & & & & & & & \\
\hline
b_2^2 & b_3^1 + b_3^3 & b_4^2 & b_4^1 + b_4^4 & & & & & & & &\\
\hline
b_2^3 &  b_3^2 + b_3^3 &  b_4^2 & b_4^2 & b_4^3 + b_4^4 & & & & & & &\\
\hline
\hline
b_3^1 & b_4^1 + b_4^2 & 0 &0 & 0& 0& & & & & &\\
\hline
b_3^2 &  b_4^2 + b_4^3 & 0 & 0& 0& 0& 0 & & & & &\\
\hline
b_3^3 & b_4^2 + b_4^4  &0 & 0& 0& 0& 0& 0& & & &\\
\hline
\hline
b_4^1 & 0 & 0 &0 & 0& 0& 0&0 & 0& & &\\
\hline
b_4^2 &  0 & 0 &0 & 0& 0& 0 &0 & 0& 0& & \\
\hline
b_4^3 & 0& 0 &0 & 0& 0& 0 &0 & 0& 0 & 0& \\
\hline
b_4^4 & 0 & 0 &0 & 0& 0& 0 &0 & 0& 0 & 0 & 0\\
\hline
\hline
b_5^1 & b_6^1  & b_7^1 & b_7^2 & 0 & b_8& 0&0 & b_9 & 0 & 0 &0\\
\hline
b_5^2 &  b_6^1 + b_6^2 + b_6^3 & b_7^2+ b_7^3& b_7^1+ b_7^3& b_7^1 + b_7^2 &b_8 & b_8&b_8  & 0  & b_9 & 0 & 0\\
\hline
b_5^3 & b_6^2 & b_7^1 & 0 & b_7^3 &0 &b_8 & 0&  0 & 0   & b_9 &0\\
\hline
b_5^4 & b_6^3 & 0 & b_7^2 & b_7^3 & 0&0 &b_8 & 0& 0 & 0   & b_9\\
\hline
\hline
b_6^1 & b_7^1 + b_7^2 & b_8 & b_8 & 0&  b_9 & 0 & 0 &0 &0 & 0&0\\
\hline
b_6^2 &  b_7^1 + b_7^3 &  b_8 & 0 & b_8   & 0  & b_9 & 0& 0& 0&0 &0\\
\hline
b_6^3 & b_7^2 + b_7^3  &  0 &  b_8  & b_8 & 0 & 0   & b_9 &0 & 0&0 &0\\
\hline
\hline
b_7^1 &  b_8 &   b_9     & 0& 0& 0&0 &0 &0 & 0& 0 &0\\
\hline
b_7^2 & b_8  & 0     & b_9& 0&0 &0 &0 &0 & 0&0 &0\\
\hline
b_7^3 &  b_8 & 0 & 0&  b_9& 0&0 &0 & 0& 0& 0& 0 \\
\hline
\hline
b_8 & b_9 & 0  &0 &0 &0 & 0&0 &0 &0 & 0&0\\
\hline
\end{array}
\end{equation*}
}

\section{The triangle inequalities for $D_4$}
 
Consider the basis
$\{\epsilon_i^*\}_{i=1,\dots, 4}$ of $\mathfrak{h}$ which is dual to the standard
basis
$\{\epsilon_i\}_{i=1,\dots, 4}$ of $\mathfrak{h}^*$ as in [Bo], Planche IV. 
Express any $h\in  \mathfrak{a}$ in this basis: 
$$h=x\epsilon_1^*+y\epsilon_2^*+z\epsilon_3^*+w\epsilon_4^*, \,\,\,x,y,z,w\in 
{\Bbb R}.$$
Then, 
$$h\in \Delta \,\,\text{iff}\,\,\, x\geq y\geq z\geq \mid w\mid.$$

\subsection{The system of inequalities corresponding to the parabolic subgroup $P_1$}

We give below the complete list (up to a permutation) of the Schubert classes
$(b_{i_1}^{j_1}, b_{i_2}^{j_2}, b_{i_3}^{j_3})=
(\epsilon^{P_1}_{w_1}, \epsilon^{P_1}_{w_2}, \epsilon^{P_1}_{w_3})$
such that
$$
\epsilon^{P_1}_{w_1}\odot_0 \epsilon^{P_1}_{w_2}\odot_0
\epsilon^{P_1}_{w_3}=\epsilon^{P_1}_{w_o^{P_1}}$$
and  write down the corresponding inequality $T^{P_1}_{{\bf w}}$:
$$\sum_{j=1}^3\,\langle \la_{w_j}^{P_1}, h_j\rangle \geq 0.$$
 We express $h_j=(x_j,y_j,z_j,w_j), j=1,2,3$ in the coordinates
$\{\epsilon_i^*\}_{i=1,\dots, 4}$.
We divide the set of inequalities in two disjoint sets, one coming from the Schubert
classes $(b_{i_1}^{j_1}, b_{i_2}^{j_2}, b_{i_3}^{j_3})$ such that at least one of the
cohomology classes is $1$. It can be seen that  the corresponding inequalities are the {\it weak triangle
inequalities} $WTI$ defined in Section 3.1.1. The remaining inequalities are called the {\it essential
triangle inequalities} $ETI$. We label the inequalities $ETI$ corresponding to
 the parabolic $P_1$ by $ETI(1)$ and similarly for $WTI$.

\bigskip 
\noindent
{\bf ETI(1):}
\begin{equation*}
\begin{array}{ccc}
(b_1, b_1, b_4):    &\quad y_1+ y_2-z_3\ge 0 &\quad (3)\\
(b_1, b_2, b_3^1):&\quad y_1+ z_2+w_3\ge 0 &\quad (6)\\
(b_2, b_3^2, b_1):&\quad z_1-w_2+y_3\ge 0 &\quad (6)
\end{array}
\end{equation*}
{\bf WTI(1):}
\begin{equation*}
\begin{array}{ccc}
(1, 1, b_6): &\quad x_1+x_2-x_3\geq 0 &\quad (3)\\
(1, b_1, b_5):&\quad x_1+y_2-y_3\geq 0 &\quad (6)\\
(1, b_2, b_4): &\quad x_1+z_2-z_3\geq 0 &\quad (6)\\
(1, b_3^1, b_3^2): &\quad x_1+w_2-w_3\geq 0 &\quad (6)\\
\end{array}
\end{equation*}
To get the full set of inequalities $T_{\bf w}^{P_1}$ for $P_1$, we need to permute the
above collection of inequalities where the subscripts $\{1,2,3\}$ are permuted
arbitrarily.
The number at the end of each inequality denotes the number of inequalities
obtained by permuting that particular inequality. Thus, the total number of inequalities
$T_{\bf w}^{P_1}$ corresponding to $P_1$ is $36$.

\subsection{The system of inequalities $T_{\bf w}^{P_2}$ corresponding to the parabolic subgroup 
$P_2$}

In each cohomological degree except for 8 and 10 there is only one orbit of the Schubert classes under
$\hat{F}$. In degree 8 there are two orbits (of three classes in the orbit of $b_4^1$ and
one in the orbit of $b_4^2$).

Of course, $\hat{F}$ acts diagonally on
the set of triples
$(b_{i_1}^{j_1}, b_{i_2}^{j_2}, b_{i_3}^{j_3})$ such that
$b_{i_1}^{j_1}\odot_0 b_{i_2}^{j_2}\odot_0b_{i_3}^{j_3}=b_9.$  Also,
$S_3$ acts on such triples via permutation and these two actions commute. So, we
get an action of the product group $S_3\times \hat{F}$ on such triples.  The following is a complete
list of such triples of Schubert classes in $H^*(G/P_2)$ up to the action of $S_3\times \hat{F}$
and the corresponding inequality $T_{\bf w}^{P_2}$. The number at the end of each inequality denotes the
order of the corresponding $S_3\times \hat{F}-$orbit.

\bigskip 

\noindent
{\bf ETI(2):}
\begin{equation*}
\begin{array}{ccc}
(b_1, b_1, b_7^1):  &\quad x_1+z_1 +    x_2+z_2     \ge y_3+z_3 &\quad (9)\\
(b_1, b_2^1, b_6^1):    &\quad x_1+z_1 +    y_2+z_2     \ge y_3+w_3 &\quad (36)\\
(b_1, b_3^1, b_5^1):    &\quad x_1+z_1 +    y_2+w_2     \ge z_3+w_3 & \quad
(18)\\
(b_1, b_3^1, b_5^2):    &\quad x_1+z_1 +    y_2+w_2     \ge y_3-z_3 &\quad (18)\\
(b_2^1, b_2^1, b_5^1):  &\quad y_1+z_1      +y_2+z_2    \ge z_3+w_3 & \quad
(18)\\
(b_2^1, b_2^2, b_5^2):  &\quad y_1+z_1      +x_2+w_2    \ge y_3-z_3  &\quad (18)\\
\end{array}
\end{equation*}

\bigskip

\noindent
{\bf WTI(2):}
\begin{equation*}
\begin{array}{ccc}
(1, 1, b_9): &\quad x_1+y_1+x_2+y_2-x_3-y_3\geq 0 &\quad (3)\\
(b_1, 1, b_8) : &\quad x_1+z_1+x_2+y_2-x_3-z_3\geq 0 &\quad (6)\\
(b_2^1, 1, b_7^1): &\quad y_1+z_1+x_2+y_2-y_3-z_3\geq 0 &\quad (18)\\
(b_3^1, 1, b_6^1): &\quad y_1+w_1+x_2+y_2-y_3-w_3\geq 0 &\quad (18)\\
(b_4^1, 1, b_5^1): &\quad z_1+w_1+x_2+y_2-z_3-w_3\geq 0 &\quad (18)\\
(b_4^2, 1, b_5^2): &\quad y_1-z_1+x_2+y_2-y_3+z_3\geq 0 &\quad (6)
\end{array}
\end{equation*}

The group  $S_3\times F$ acts canonically on $\mathfrak{a}^3$, where $S_3$ acts by
permutation of the three factors
and $F$ acts via the corresponding Dynkin automorphism of $\mathfrak{a}$.
To get the full set of inequalities $T_{\bf w}^{P_2}$ for $P_2$, we need to apply the group
 $S_3\times F$ to the above collection of inequalities.

Thus, we get totally $186$ inequalities corresponding to the maximal parabolic $P_2$.

\vskip1ex

The multiplication table for $H^*(G/P_3)$ (resp. $H^*(G/P_4)$) can be obtained from that of
$H^*(G/P_1)$ by using the isomorphism of $H^*(G/P_1)$ with $H^*(G/P_3)$
(resp. $H^*(G/P_4)$) induced from  the Dynkin automorphisms. Accordingly, the 
inequalities corresponding to $H^*(G/P_3)$ and $H^*(G/P_4)$ are obtained from  
$T_{\bf w}^{P_1}$ by applying the action of $F$. All in all, each system 
$T_{\bf w}^{P_3}$ and $T_{\bf w}^{P_4}$ consists of 
36 inequalities.

Below are the explicit lists of inequalities.

\subsection{The system of inequalities corresponding to the parabolic subgroup $P_3$}

\bigskip ~

{\bf ETI(3):}
\begin{equation*}
\begin{array}{c}
x_1+y_1-z_1+w_1+x_2+ y_2-z_2
+w_2-x_3+y_3-z_3-w_3\ge 0 \\
x_1+y_1-z_1+w_1+x_2-y_2+z_2
+w_2-x_3+y_3+z_3+w_3\ge 0\\
x_1-y_1+z_1+w_1+x_2- y_2-z_2
-w_2+x_3+y_3-z_3+w_3\ge 0 
\end{array}
\end{equation*}

{\bf WTI(3):}
\begin{equation*}
\begin{array}{cc}
x_1+y_1+z_1-w_1+x_2+ y_2+z_2
-w_2-x_3-y_3-z_3+w_3\ge 0  \\
x_1+y_1+z_1-w_1+x_2+ y_2-z_2
+w_2-x_3-y_3+z_3-w_3\ge 0 \\
x_1+y_1+z_1-w_1+x_2- y_2+z_2
+w_2-x_3+y_3-z_3-w_3\ge 0 \\
x_1+y_1+z_1-w_1-x_2+ y_2+z_2
+w_2+x_3-y_3-z_3-w_3\ge 0\\
\end{array}
\end{equation*}
To get the full set of inequalities $T_{\bf w}^{P_3}$ for $P_3$, we need to permute the
above collection of inequalities where the subscripts $\{1,2,3\}$ are permuted
arbitrarily.

\subsection{The system of inequalities corresponding to the parabolic subgroup $P_4$}

\bigskip
~

\noindent {\bf ETI(4):}
\begin{equation*}
\begin{array}{cc}
x_1+y_1-z_1-w_1+x_2+ y_2-z_2
-w_2-x_3+y_3-z_3+w_3\ge 0 \\
x_1+y_1-z_1-w_1+x_2-y_2+z_2
-w_2+x_3-y_3-z_3+w_3\ge 0 \\
x_1-y_1+z_1-w_1-x_2+y_2+z_2
-w_2+x_3+y_3-z_3-w_3\ge 0 
\end{array}
\end{equation*}
{\bf WTI(4):}
\begin{equation*}
\begin{array}{cc}
x_1+y_1+z_1+w_1+x_2+ y_2+z_2
+w_2-x_3-y_3-z_3-w_3\ge 0 \\
x_1+y_1+z_1+w_1+x_2+ y_2-z_2
-w_2-x_3-y_3+z_3+w_3\ge 0 \\
x_1+y_1+z_1+w_1+x_2- y_2+z_2
-w_2-x_3+y_3-z_3+w_3\ge 0 \\
x_1+y_1+z_1+w_1+x_2- y_2-z_2
+w_2-x_3+y_3+z_3-w_3\ge 0\\
\end{array}
\end{equation*}
To get the full set of inequalities $T_{\bf w}^{P_4}$ for $P_4$, we need to permute the
above collection of inequalities where the subscripts $\{1,2,3\}$ are permuted
arbitrarily.

\subsection{The cone ${\mathcal C}$}

Thus, the total number of inequalities $T_{{\bf w}}^{P_i}$ defining
the cone ${\mathcal C}$ inside $\Delta^3$ is equal to $36+186+36+36
=294$.  Since $\Delta^3 \subset
\mathfrak{a}^3$ is defined by $12$ inequalities, we get altogether
$306$ inequalities defining the cone ${\mathcal C}$ inside
$\mathfrak{a}^3$. Let $\Sigma$ be the set of these  $306$
inequalities defining the cone ${\mathcal C}$.

\begin{thm}
The system $\Sigma$ is irredundant.
\end{thm}
\proof In order to show the irredundancy of the system $\Sigma$, it suffices to show that the cone
  ${\mathcal C}$
has 306 faces.
It is done by applying the MAPLE package CONVEX  \cite{convex} to the above  system (see
\cite{kapovich1}). \qed

\begin{rem}
The same computation also shows that the cone ${\mathcal C}$  has
81 extremal rays.
\end{rem}

Our next goal is to verify the {\em saturation conjecture} for the
group $G^\vee=$ Spin$(8)$. Let $P^\vee\subset ({\mathfrak h}^\vee)^*={\mathfrak h}$ denote the
weight lattice of $G^\vee$ and $Q^\vee\subset P^\vee$ denote the root lattice. Of course,
$Q=Q^\vee$ since $G$ is simply-laced.
Recall that in the Introduction we have defined  the semigroups
$\mathcal{C}^0$ and $\mathcal{R}$ of ${\mathcal D}^3$ with
$$\mathcal{R} \subset \mathcal{C}^0.$$

\begin{thm}
[Saturation theorem for $Spin(8)$]
\label{saturation}
$$
\mathcal{R}= \mathcal{C}^0. $$
\end{thm}
\proof In order to prove the inclusion $\mathcal{C}^0 \subset \mathcal{R}$, it suffices to show that
each semigroup generator of $\mathcal{C}^0$ belongs to $\mathcal{R}$.  To find the
minimal set of semigroup generators (Hilbert basis)
for $\mathcal{C}^0$ , we define a basis $\{\bar{\al_i}, \zeta_j\}_{1\leq i\leq 4,1\leq j\leq 8}$
 of the lattice
$
\phi^{-1}(Q^\vee),
$
where
$$
\phi: (P^\vee)^3\to P^\vee, \,\,\phi(\la,\mu,\nu)=\la+\mu+\nu.
$$

Consider the splitting of the exact sequence (for $K:=\text{Ker}\, \phi$)
$$
0\to K\to (P^\vee)^3 \stackrel{\phi}{\to} P^\vee \to 0
$$
over $Q^\vee$ under the map
$
\psi (\al_i)=(\al_i, 0, 0), i=1,...,4.$

Therefore, we can identify $Q^\vee$ with the subgroup $\psi(Q^\vee)\subset (P^\vee)^3$ with basis
$$
\{\bar{\alpha_i}=\psi(\alpha_i)\}_{1\leq i \leq 4}. 
$$
Set 
$$
\zeta_j=(-\omega_j, \omega_j,0)$$
 for $1\leq j\leq 4$ and
$$
\zeta_j=(-\omega_{j-4}, 0,\omega_{j-4})$$
 for $5\leq j\leq 8$, where $\{\omega_j\}$
are the fundamental weights for
$Spin (8)$. Then, it is clear that
$
\{\zeta_j, j=1,...,8\}$ is a basis of $K$ and
$$
\{\bar{\al_i}, \zeta_j : i=1,...,4, j=1,...,8 \}$$
is a basis of $\phi^{-1}(Q^\vee)$.

Thus, the semigroup $\mathcal{C}^0$  is precisely equal to the integral points of the cone
 ${\mathcal C}$  with respect to the coordinates  $\{\bar{\al_i}, \zeta_j\}$.

Computation of the Hilbert basis $H$ of $\mathcal{C}^0$  is done
via the package HILBERT, \cite{42} (see \cite{kapovich2}). Observe
that the action of the group $ S_3\times F $ on $\mathfrak{a}^3$
keeps $\mathcal{C}^0$ stable. Since the Hilbert basis is unique,
it follows that $H$  is invariant under the action of $S_3\times F$.
Below is the list $H'$ of elements of $H$ modulo the action of
$S_3\times F$:
$$
\begin{array}{c}
(\omega_1, \omega_1, 0)\\
(\omega_2, \omega_2, 0)\\
(\omega_2, \omega_2, \omega_2)\\
(\omega_1, \omega_3, \omega_4)\\
(\omega_1, \omega_1, \omega_2)\\
(\omega_1, \omega_2, \omega_3+\omega_4)\\
(2\omega_1, \omega_2, \omega_2)\\
(\omega_1+\omega_2, \omega_2, \omega_3+\omega_4)\\
(\omega_2, \omega_2, \omega_1+\omega_3+\omega_4)\\
(2\omega_2, \omega_2, \omega_1+\omega_3+\omega_4).
\end{array}
$$

Since $S_3\times F$ also preserves the semigroup $\mathcal{R}$, in
order to prove Theorem \ref{saturation}, it suffices to check that
$H'\subset \mathcal{R}$. This is done
 using  MAPLE package WEYL written by John Stembridge, see \cite{stembridge}.
It is done in \cite{kapovich2}.


Addresses:

M.K.: Department of Mathematics, University of California,
Davis, CA 95616, USA. (kapovich@math.ucdavis.edu)

S.K.: Department of Mathematics, University of North Carolina,
Chapel Hill, NC 27599-3250, USA. (shrawan@email.unc.edu)

J.M.: Department of Mathematics, University of Maryland, College
Park, MD 20742, USA. (jjm@math.umd.edu)
\end{document}